\documentclass[12pt]{article}
\usepackage{graphicx}
\usepackage{amssymb}
\usepackage{amsmath}
\usepackage{float}
\newcommand{\eps}{\varepsilon}

\numberwithin{equation}{section}
\def\la{{\lambda}}
\def\BC{{\mathbb C}}
\def\BR{{\mathbb R}}

\newtheorem{Pa}{Paper}[section]
\newtheorem{Tm}[Pa]{{\bf Theorem}}
\newtheorem{La}[Pa]{{\bf Lemma}}
\newtheorem{Cy}[Pa]{{\bf Corollary}}
\newtheorem{Rk}[Pa]{{\bf Remark}}

\newtheorem{Ee}[Pa]{{\bf Example}}
\newtheorem{Dn}[Pa]{{\bf Definition}}
\newtheorem{Pn}[Pa]{{\bf Proposition}}

\def\XXint#1#2#3{{\setbox0=\hbox{$#1{#2#3}{\int}$}
     \vcenter{\hbox{$#2#3$}}\kern-.5\wd0}}

\title{Almost periodic functions and \\ an analytical method of solving \\
the number partitioning problem}

\author{Lev Sakhnovich}

\date{}

\begin{document}
\maketitle

\begin{abstract} In the present  paper, we study the limit sets of the almost periodic functions $f(x)$.
It is interesting  that the values $r=\inf|f(x)|$ and $R=\sup|f(x)|$ may be expressed in the exact form.
We show that the ring $r\leq |z|\leq R$ is the limit set of the almost periodic function $f(x)$ (under some natural conditions on $f$).
The exact expression for $r$ coincides with the well known partition problem formula and gives a new analytical method  of solving the
corresponding partition problem.
Several interesting examples are considered.
For instance, in the case of the five numbers,
the well-known (see, e.g., \cite{BoMer}) Karmarkar--Karp algorithm  gives the value $m=2$ as the solution of the
partition problem in our example,  and our method gives the correct answer $m=0.$
The figures presented in Appendix illustrate  our results. 
\end{abstract}



\section{Introduction}

In the present  paper, we study the limit sets $\Delta_{\pm}$ of the almost periodic functions
\begin{equation}f(x)=\sum_{k=1}^{n}c_{k}e^{ix\lambda_k} \quad (\la_k\in \BR,\,\, c_k\in \BC),\label{1.1}
\end{equation} 
where $\BR$ is the real axis and $\BC$ is the complex plane.

We proved the following assertion:

If the real numbers $\lambda_1,\lambda_2,...\lambda_n$ are linearly independent over the field of rational numbers, then $\Delta_{\pm}$ coincides with the ring $r\leq|z|\leq{R}$, where $r=\inf|f(x)|$ and $R=\sup|f(x)|$.

It is interesting, that $r$ and $R$ may be expressed in the exact form  in terms of the absolute values $|c_k|$:
\begin{equation}R=\sum_{k=1}^{n}|c_{k}|,\quad r=\min\left|\sum_{k=1}^{n}(|c_k|\varepsilon_k)\right| \quad (\varepsilon_k=\pm{1}),\label{1.2}\end{equation}
for $f$ of the form \eqref{1.1}.
The formula \eqref{1.2}  for $r$ coincides with the well known partition problem formula (compare \eqref{1.2} and \eqref{2.45}).
The equality $r=\inf|f(x)|$ presents a new analytical method for solving this  partition problem.
In particular, we consider  Example \ref{Example 3.9} of the partition problem where $n=5$ in \eqref{2.45}. The 
solution of the partition problem in Example~\ref{Example 3.9} is $m=0.$
The well-known Karmarkar--Karp algorithm (see \cite{BoMer}) gives the answer $m=2$, and our method gives the correct
solution.
In the last section~4, we consider the  case where
\begin{equation}f(x)=\sum_{k=1}^{\infty}c_{k}e^{ix\lambda_k},\label{1.3}\end{equation}
and study the sequence $m_{n}=\min|\sum_{k=1}^{n}(|c_k|\varepsilon_k)|$ for $n{\to}\infty$.
The figures presented in Appendix illustrate  our results.
\section{Definitions and properties of almost periodic functions}
\setcounter{equation}{0}
1. Now, we introduce the main notions and properties of the almost periodic functions \cite{Bes, Bor, LZ}.
\begin{Dn}\label{Definition 2.1}The set $E\subseteq \BR$ is called relatively dense  if there exists a number $L>0$ such that any interval of length $L$ contains at least one number from
$E$.\end{Dn}
\begin{Dn}\label{Definition 2.2} Let $f(x)$ be  defined on $\BR$. The number $\tau$ is called
a translation number of $f(x)$ belonging to the given $\varepsilon{\geq}0$ if
\begin{equation}|f(x+\tau)-f(x)|{\leq}\varepsilon \quad (-\infty<x<\infty).\label{2.1}\end{equation}\end{Dn}
We denote the set of all the translation numbers of a function $f(x)$ belonging to $\varepsilon$ by $E\big(\varepsilon,f(x)\big)$.
\begin{Dn}\label{Definition 2.3}A continuous function is uniformly almost periodic if for any $\varepsilon$
the set $E\big(\varepsilon,f(x)\big)$  is relatively dense. \end{Dn}
\begin{Dn}\label {Definition 2.4}A finite set of real numbers $\lambda_1,\lambda_2...\lambda_n$ is called linearly independent over the field of rational numbers if the equality
\begin{equation}\sum_{k=1}^{n}r_k\lambda_k=0, \label{2.2}\end{equation}
 for rational $r_{k}$, implies that all $r_k$ are equal to zero.\end{Dn}
\begin{Tm}(Kronecker's theorem) Let two sets of real numbers $\lambda_1,\lambda_2, \ldots  , \lambda_n$ and $\theta_1,\theta_2, \ldots , \theta_n$ be given.  Then, the inequalities
\begin{equation}|\lambda_k\tau-\theta_k|<\delta \quad (\mod \,2\pi)\quad (1\leq k\leq n)\label{2.3}\end{equation}
have a consistent  real solution $\tau$ for any positive number $\delta$ if and only if the numbers
$\lambda_1,\lambda_2,...\lambda_n$ are linearly independent over the field of rational numbers.\end{Tm}
2. Further, we investigate  the almost periodic function of the form
\begin{equation}f(x)=\sum_{k=1}^{n}c_{k}e^{ix\lambda_k}.\label{2.4}\end{equation}

Let us introduce the limit sets  $\Delta_{+}$ and $\Delta_{-}$ for the function $f(x)$ in a  precise way.
Namely, the number $\alpha$ belongs to $\Delta_{+}$ if there is a sequence ${x_k}$ such that
$$\alpha=\lim_{k\to\infty}f(x_k), \quad x_k{\to}+\infty.$$
The number $\alpha$ belongs to $\Delta_{-}$ if there is a sequence ${x_k}$ such that
$$\alpha=\lim_{k\to\infty}f(x_k), \quad x_k{\to}-\infty.$$
The set of all limit points is denoted by 
$\Delta:=\Delta_{+}{\cup}\Delta_{-}$.
The following theorem is valid for the limit sets of the almost periodic functions  $f$ of the form \eqref{2.4}.
\begin{Tm}\label{Theorem 2.6}Let the set  of real numbers $\lambda_1,\lambda_2,\ldots ,\lambda_n$ be linearly independent over the field of rational numbers and let $f(x)$ have
the form \eqref{2.4}. Then, each point $z$ on the circle $|z|=|f(x)|$ belongs to the limit  sets $\Delta_{-}=\Delta_{+}=\Delta$ for $f$.\end{Tm}
\emph{Proof.} We assume that all $\theta_k$ are equal to $\theta.$ Then we have
\begin{equation}|e^{i\theta}f(x)-f(x+\tau)|{\leq}\sum_{k=1}^{n}|c_{k}||e^{i\theta}-e^{i\tau\lambda_k}|.\label{2.5}\end{equation}
If $\delta_N{\to}0$  then $f(x+\tau_N){\to}e^{i\theta}f(x)$ (see \eqref{2.3}-\eqref{2.5}). The theorem is proved.\\

Let us introduce the following notations:
\begin{equation}r=\inf|f(x)|,\quad R=\sup|f(x)| \quad (-\infty<x<\infty).\label{2.6}\end{equation}
Theorem \ref{Theorem 2.6} implies the assertion:
 \begin{Cy}\label{Corollary 2.7} Let the conditions of Theorem \ref{Theorem 2.6} be fulfilled. Then,
 the limit set $\Delta$ of the function $f(x)$ coincides with the ring
 \begin{equation}r{\leq}|z|{\leq}R.\label{2.7}\end{equation}\end{Cy}
 \begin{Ee}\label{Example 2.8} Consider the case 
 \begin{equation}f(x)=e^{ix}+e^{i\sqrt{2}x}, \quad -\infty<x<\infty.\label{2.8}\end{equation}\end{Ee}
 The function $f(x)$ can be represented in the  parametric form
 \begin{equation}x=\cos{t}+\cos\sqrt{2}\, t,\quad y=\sin{t}+\sin{\sqrt{2}\, t}. \label{2.8+}\end{equation}
  The corresponding parametric plot is illustrated by Figure 1 in Appendix.
We see that
 \begin{equation}r=0,\quad R=2.\label{2.10}\end{equation}
  \begin{Ee}\label{Example 2.9}. Consider the case 
 \begin{equation}f(x)=e^{ix}+(1/2)e^{i\sqrt{2}x}, \quad -\infty<x<\infty.\label{2.11}\end{equation}\end{Ee}
 The function $f(x)$ can be represented in the  parametric form
\begin{equation}x=\cos{t}+(1/2)\cos\sqrt{2}\, t,\quad y=\sin{t}+(1/2)\sin{\sqrt{2}\, t}.
 \label{2.12}\end{equation}
  The corresponding parametric plot is illustrated by Figure 2 in Appendix.
We see that
\begin{equation}r=1/2,\quad R=3/2.\label{2.13}\end{equation}

Note that $\lambda_1=1$ and $\lambda_2=\sqrt{2}$, which are considered in   Examples \ref{Example 2.8} and \ref{Example 2.9}, are linearly 
independent over the field of rational numbers.

3. \emph{Now, we will  study the behavior of the radii $r$ and $R$.}
\begin{Tm}\label{Theorem 2.10} Let the set of real numbers $\lambda_1,\lambda_2, \ldots , \lambda_n$ be linearly independent over the field of rational numbers and 
let $f(x)$ have the form \eqref{2.4}. Then,
\begin{equation}R=\sum_{k=1}^{n}|c_k|.\label{2.14}\end{equation} \end{Tm}
\emph{Proof.} We assume that all $\theta_k=-\arg(c_k)$. Then we have
\begin{equation}|g(x)-f(x+\tau)|{\leq}\sum_{k=1}^{n}|c_{k}||e^{i\theta_k}-e^{i\tau\lambda_k}|,\label{2.15}\end{equation}
where
\begin{equation}g(x)=\sum_{k=1}^{n}|c_{k}|e^{ix\lambda_k}.\label{2.16}\end{equation}

If $\delta_N{\to}0$  then $f(x+\tau_N){\to}g(x)$ (see \eqref{2.3},\eqref{2.15}, \eqref{2.16}).
Hence
\begin{equation}\sup|f(x)|=\sup|g(x)|=g(0)=\sum_{k=1}^{n}|c_k|.\label{2.17}\end{equation}
\begin{equation}\inf|f(x)|=\inf|g(x)|.
\label{2.18}\end{equation}
The theorem is proved.\\

\begin{La}\label{Lemma 2.10}
Let the set of real numbers $\lambda_1,\lambda_2, \ldots , \lambda_n$ be linearly independent over the field of rational numbers and  let $f(x)$ have the form \eqref{2.4}.  Then,
\begin{equation}r{\leq}\min\Big|\sum_{k=1}^{n}(|c_k|\varepsilon_k)\Big|,\quad \varepsilon_k=\pm{1}.\label{2.19}\end{equation}\end{La}
\emph{Proof.}  The values of the radii r and R depend only on the modules $|c_k|$ (see \eqref{2.17} and \eqref{2.18}).
Hence the function
\begin{equation}G(x)=\sum_{k=1}^{n}|c_{k}|{\varepsilon_k}e^{ix\lambda_k},\quad \varepsilon_k=\pm1\label{2.20}\end{equation} has the same r as the function $f(x)$. We have
\begin{equation}G(0)=\sum_{k=1}^{n}|c_{k}|{\varepsilon_k},\quad \varepsilon_k=\pm1\label{2.21}\end{equation}
Relation \eqref{2.19} follows from equality \eqref{2.21}. Lemma is proved.
\begin{La}\label{Lemma 2.12}Let the set of real numbers $\lambda_1,\lambda_2,\dots\lambda_n$ be linearly independent over the field of rational numbers,
let $f(x)$ have the form \eqref{2.4}, and assume that $r>0$. Then,
\begin{equation}r{\geq}\min\Big|\sum_{k=1}^{n}(|c_k|\varepsilon_k)\Big|,\quad \varepsilon_k=\pm{1}.\label{2.22}\end{equation}\end{La}
\emph{Proof.} We shall consider the case where $\varepsilon_k=1.$ In this case, $G(x)=G(x,\lambda)$ has the form
\begin{equation}G(x,\lambda)=\sum_{k=1}^{n}|c_{k}|e^{ix\lambda_k},\label{2.23}\end{equation}
where $\lambda=(\lambda_1,\lambda_2,\ldots,\lambda_n)$.
 Here, $\lambda=(\lambda_1,\lambda_2,\ldots,\lambda_n)$ and $x$ are fixed.
 
It follows from the relation \eqref{2.23} that
\begin{equation} |G(x,\lambda)|^{2}=A^{2}(x,\lambda)+B^{2}(x,\lambda), \label{2.24}\end{equation}
where
\begin{equation} A(x,\lambda)=\sum_{k=1}^{n}|c_{k}|\cos(x\lambda_k),\quad B(x,\lambda)=\sum_{k=1}^{n}|c_{k}|\sin(x\lambda_k).\label{2.25}\end{equation}
 To find the extremal points of $|G(x,\lambda)|$ we consider the relations
\begin{equation}\frac{\partial}{\partial{\lambda_j}}|G(x,\lambda)|^{2}=-2xA(x,\lambda)|c_{j}|\sin(x\lambda_j)+
2xB(x,\lambda)|c_{j}|\cos(x\lambda_j)=0.\label{2.26}\end{equation}
According to our assumption, the functions $A (x,\lambda)$ and $B (x,\lambda)$ cannot be equal to zero simultaneously. We begin with the case where 
$$A(x,\lambda)=0,\quad B(x,\lambda){\ne}0.$$
We have
$\cos(x\lambda_j)=0.$ So, in this case,  the stationary points of the function $G(x,\la)$ are defined by the relations
$x\lambda_j=\pi/2+n_j\pi$,
 where $n_j$ are integer numbers. In the stationary points, the function $G(x,\lambda)$ takes the values
\begin{equation} G(x,\lambda_{st})=i\sum_{k=1}^{n}|c_{k}|\varepsilon_k,\quad \varepsilon_k=(-1)^{n_k}.
\label{2.27}\end{equation}
Now, we assume that $B(x,,\lambda)=0,\,A(x,\lambda){\ne}0.$
We have
$\sin(x\lambda_j)=0.$ So, in this case  the stationary points of the function $G(x,,\lambda)$ are defined by the relations
$x\lambda_j=n_j\pi$,
 where $n_j$ are integer numbers. In the stationary points, the function $G(x,\lambda)$ takes the values
\begin{equation} G(x,\lambda_{st})=\sum_{k=1}^{n}|c_{k}|\varepsilon_k,\quad \varepsilon_k=(-1)^{n_k}.
\label{2.28}\end{equation}
The last case: $A(x,\lambda){\ne}0,\,B(x,\lambda){\ne}0.$ It follows from \eqref{2.26} that
\begin{equation}\tan(x\lambda_j)=B(x,\lambda)/A(x,\lambda).\label{2.29}\end{equation}
Hence, we have
\begin{equation}x\lambda_j=\alpha(x)+n_{j}\pi,\label{2.30}\end{equation}
where $n_j$ are integer numbers. Thus, in view of \eqref{2.23}, the function $G(x,\lambda)$  satisfies in the stationary points the relation 
\begin{equation} G(x,\lambda_{st})=e^{i\alpha(x)}\sum_{k=1}^{n}|c_{k}|\varepsilon_k,\quad \varepsilon_k=(-1)^{n_k}.
\label{2.31}\end{equation} We found all the stationary points. All extremal points are stationary points.
Hence, it is proved that
\begin{equation}r=\inf|G(x,\lambda)|{\geq}\inf|G(x,\lambda_{st})|=\min\Big|\sum_{k=1}^{n}|c_{k}|\varepsilon_k\Big|.
\nonumber\end{equation}where $-\infty<x<+\infty.$
 We obtained the assertion of the lemma.
\begin{Cy}\label{Corollary 2.13}Let $f(x)$ have the form \eqref{2.4} and set
\begin{equation}M=\sup|f(x)|,\quad m=\inf|f(x)|\quad (-\infty<x<+\infty).\label{2.32}\end{equation}
Then,
\begin{equation}M{\leq}\sum_{k=1}^{n}|c_{k}|,\quad m{\geq}\min|\sum_{k=1}^{n}(|c_k|\varepsilon_k)|,\quad \varepsilon_k=\pm{1}.\label{2.33}\end{equation}\end{Cy}
Let us return to the case where the real numbers $\lambda_1,\lambda_2,\ldots, \lambda_n$ are linearly independent.
Relations \eqref{2.14}, \eqref{2.19} and \eqref{2.22} imply the assertion:
\begin{Tm}\label{Theorem 2.14}
Let the set of the real numbers $\lambda_1,\lambda_2,\ldots,\lambda_n$ be linearly independent over the field of rational numbers and let $f(x)$ have the form \eqref{2.4}. Then,
\begin{equation}M=R=\sum_{k=1}^{n}|c_{k}|,\quad m=r=\min|\sum_{k=1}^{n}(|c_k|\varepsilon_k)|,\quad \varepsilon_k=\pm{1}.\label{2.34}\end{equation}\end{Tm}
\emph{Proof.} The first equality in \eqref{2.34} was proved in Theorem \ref{Theorem 2.10}.  The second equality in \eqref{2.34} follows (for $r>0$) from Lemmas \ref{Lemma 2.10} and \ref{Lemma 2.12}. Hence, this relation is valid for
$r=0$ as well. The theorem is proved.\\

In the qualitative theory of differential equations, the limit sets of trajectories are studied  in the proof of Poincare-Bendixson theorem (see \cite{Lef}). Applying these arguments to our case, we obtain the assertion:
\begin{Pn}\label{Proposition 2.15}Let the function $f(x)$ be of the form \eqref{2.4}. Then, the corresponding sets $\Delta_{+}$ and  $\Delta_{-}$ are closed, connected and nonempty.\end{Pn}
4. \emph{Periodical case.}
Assume that $f(x)$ has the form:
\begin{equation}f(x)=\sum_{k=1}^{n}c_{k}e^{ixk\alpha},\quad n{\geq}2.\label{2.35}\end{equation}
The numbers $\lambda_k=k\alpha$ are linearly dependent and the corresponding function $f(x)$ is periodic:
$f(x)=f(x+2\pi/{\alpha})$.
\begin{Ee}\label{Example 2.16} Consider the case where
\begin{equation}f(x)=e^{ix} +ae^{2ix}, \quad a>0.\label{2.36}\end{equation}\end{Ee}
Let us set
\begin{equation}A(x)=\cos(x)+a\cos(2x),\quad B(x)=\sin(x)+a\sin(2x).\label{2.37}\end{equation}
Hence, the equality
 \begin{equation}|G(x)|^2=A^2(x)+B^2(x)=1+a^2+2a\cos(x)\label{2.38}\end{equation}
holds. Relation \eqref{2.38} implies that
\begin{equation}M=1+a,\quad m=|1-a|.\label{2.39}\end{equation}
The function $f(x)$ may be represented in the  parametric form
\begin{equation}x=\cos{t}+a\cos{2t},\quad y=\sin{t}+a\sin{2t}.
 \nonumber\end{equation}
The corresponding parametric plot is illustrated, for $a=2$, in Appendix (see Figure 3).
We see that
\begin{equation}r=1,\quad R=3.\nonumber\end{equation}
5. Equality \eqref{2.39} suggests that the following fact holds:
\begin{Tm}\label{Theorem 2.17}Let the function $f(x)$ have the form \eqref{2.4}. Then,
\begin{equation}M=\sum_{k=1}^{n}|c_{k}|,\quad m=\min|\sum_{k=1}^{n}(|c_k|\varepsilon_k)|,\quad \varepsilon_k=\pm{1}.\label{2.40}\end{equation}\end{Tm}
\emph{Proof.} We introduce the functions
\begin{equation}f_{\ell}(x)=\sum_{k=1}^{n}c_{k}e^{i\lambda_{k}(\ell)x},\label{2.41}\end{equation}
where the  real numbers $\lambda_1(\ell),\lambda_2(\ell),\ldots,\lambda_n(\ell)$
are linearly independent over the field of rational numbers and $\lambda_k(\ell){\to}\lambda_k$
when $\ell{\to}\infty$. Hence we have:
\begin{equation}\sup|f_{\ell}(x)|{\to}\sup|f(x)|,\quad \inf|f_{\ell}(x)|{\to}\inf |f(x)|,\quad \ell{\to}\infty.
\label{2.42}\end{equation}
The values $\sup|f_{\ell}(x)|$ and  $\inf|f_{\ell}(x)|$ are the constants $M_{\ell}$ and $m_{\ell}$
(see Theorem \ref{Theorem 2.14}).The assertion of the theorem follows from \eqref{2.42}.
\begin{Cy}\label{Corollary 2.18} The equalities  \eqref{2.40} are valid, in particular, for a periodic $f(x)$.\end{Cy}
\begin{Ee}\label{Example 2.19} Consider the case where $f(x)$ has the form
\begin{equation}f(x)=e^{ix}+e^{2ix}+(1/10)e^{i\sqrt{3}x}.\label{2.43}\end{equation}\end{Ee}
The example under consideration is  intermediate in the following sense: the function $f(x)$ is not periodic, but the system $\lambda_1=1,\,\lambda_2=2,\,\lambda_3=\sqrt{3}$
is not linearly independent. Taking into account \eqref{2.40}, we have:
\begin{equation}M=2.1,\quad m=0.1.\label{2.44}\end{equation}
The function $f(x)$ may be represented in the  parametric form
\begin{equation}x=\cos{t}+\cos{2}t+(1/10)\cos\sqrt{3}t,\quad y=\sin{t}+\sin{2}t+(1/10)\sin\sqrt{3}t.
 \label{2.12+}\end{equation}
The corresponding parametric plot is illustrated by Figure 4.

\section{Partition problem}
\setcounter{equation}{0}
\subsection{Problem and examples}
Let a number set $\{c_1,c_2,...,c_n\}$ be given.  The  partition problem is the problem of finding the following minimum:
\begin{equation} m=\min\Big|\sum_{k=1}^{n}(|c_k|\varepsilon_k)\Big|,\quad \varepsilon_k=\pm{1},\label{2.45}\end{equation}
where the minimum is taken among the values obtained for various choices of $\{\eps_k\}$.
It is known that the  partition problem is important in number theory, in theoretical computer science and in statistical physics
\cite{BoMer}. In the present paper, we showed that the   partition problem plays  an essential role in the theory of almost periodic functions (see section 2). Further, we give some
useful examples, where the partition problem has explicit solutions, and compare our approach with the well known Karmarkar--Karp algorithm.

\emph{The following procedure is used in Karmarkar--Karp method of solving partition problem \cite{BoMer}.
A list of $n$ positive numbers is considered. The two largest ones are replaced by their difference, which yields a new list of $n-1$ numbers. Iterating this operation $n-1$ times we obtain a single number - either approximate or precise solution of the problem.}
\begin{Ee}\label{Example 3.1} Consider the set  $\{|c_1|,|c_2|\}$ $(n=2)$.\\
Then, we have  $m=\big|[|c_1|-|c_2|]\big|.$ \end{Ee}
\begin{Ee}\label{Example 3.2}  Consider the set  $\{|c_1|,|c_2|, |c_3|\}$ $(n=3)$ where $|c_1|{\geq}|c_2|{\geq}|c_3|$.\end{Ee}
Clearly, the corresponding partition problem \eqref{2.45} may be reduced to one of the two partition problems for $n=2$.
Either it is the partition problem for the set $\{|c_1|-|c_2|,|c_3|\}$ or for the set $\{|c_1|+|c_2|,|c_3|\}$.
Using Example 3.1 we obtain
$m=\big|[|c_1|-|c_2|-|c_3|]\big|$ in the first case and $m_1=\big|[|c_1|+|c_2|-|c_3|]\big|$ in the second case.
It is easy to see that $m_1{\geq}m.$ Hence, we proved that
\begin{equation}m=\min\Big|\sum_{k=1}^{3}(|c_k|\varepsilon_k)\Big|=\big|[|c_1|-|c_2|-|c_3|]\big|.
\label{3.2}\end{equation}
\begin{Rk}\label{Remark 3.3} Thus, in the case of three positive numbers, we take the largest number and subtract two others.
The absolute value  of the obtained expression equals $m$. \end{Rk}
\begin{Ee}\label{Example 3.4} Consider the set $\{|c_1|,|c_2|, |c_3|,|c_4|\}$ $(n=4)$, \\ 
where $|c_1|{\geq}|c_2|{\geq}|c_3|{\geq}|c_4|.$ \end{Ee}
Similar to the Example \ref{Example 3.2}, the corresponding partition problem may be reduced to one of the
two partition problems for  $n-1=3$ numbers, that is, 
1) for the set $\{|c_1|-|c_2|,|c_3|, |c_4|\}$ or 2) for the set $\{|c_1|+|c_2|,|c_3|,|c_4|\}.$

In the first case, using Remark \ref{Remark 3.3}, we obtain
$m=\big|[(|c_1|-|c_2|)-|c_3|-|c_4|]\big|$ if $|c_1|-|c_2|{\geq}|c_3|$ and 
$m=\big|[|c_3|-(|c_1|-|c_2|)-|c_4|]\big|$ if $|c_1|-|c_2|{\leq}|c_3|$.

In the second case,  we have
$m_1=\big|[|c_1|+|c_2|-|c_3|-|c_4|]\big|$
It is easy to see that $m_1{\geq}m.$ Hence, the following proposition is  proved for Example \ref{Example 3.4}.
\begin{Pn}\label{Proposition 3.4}If $|c_1|-|c_2|{\geq}|c_3|$, then  $m=\big| |c_1|-|c_2|-|c_3|-|c_4|\big|$.\\
If $|c_1|-|c_2|{\leq}|c_3|$, then $m=\big| |c_3|-|c_1|+|c_2|-|c_4|\big|$.\end{Pn}

 The precise solutions $m$ of the partition problems for the cases $n=2,3,4$ were  found  using  Karmarkar--Karp method.
\begin{Rk}\label{Remark3.5} For the cases $n=2,3,4,$  the Karmarkar--Karp method gives the precise results. \end{Rk}

\begin{Ee}\label{Example 3.6}  Let us consider the case $n=5$ and the set
$\{5,5,6,7,9\}$ \end{Ee} In this case, the Karmarkar--Karp method gives an approximate  result $m=2$.
The  precise result for the solution of the partition problem  is $m=0$.
(See the details in Example \ref{Example 3.9}.)
\begin{Rk}\label{Remark 3.6} For $n{\geq}5$, the Karmarkar--Karp method gives the approximate results. \end{Rk}.

\subsection{Analytic solution}
Recall the equality \eqref{2.24} for $G(x)=G(x,\la)$ given by \eqref{2.23}, $A(x)=A(x,\la)$ and  $B(x)=B(x,\la)$:
\begin{equation} |G(x)|^{2}=A^{2}(x)+B^{2}(x). \label{3.4}\end{equation}
According to Theorem \ref{Theorem 2.17}, we have
\begin{equation}
 \inf{|G(x)|}=\min\Big|\sum_{k=1}^{n}(|c_k|\varepsilon_k)\Big|,\quad -\infty<x<+\infty,\quad \varepsilon_k=\pm{1}.\label{3.5}\end{equation}
Thus, the problem of discrete analysis (Partition problem) is reduced to the problem of continuous analysis
($\inf|G(x)|$).
\begin{Ee}\label{Example 3.9} We return to the set  $\{5,5,6,7,9\}$ $(n=5)$ mentioned in Example \ref{Example 3.6}.
\end{Ee}
Let us introduce the function
\begin{equation} f(x,a)=5e^{ixa}+5e^{ixa\sqrt{2}}+6e^{ixa\sqrt{3}}+7e^{ixa\sqrt{5}}+9e^{ixa\sqrt{7}}.\label{3.4+}
\end{equation}
The calculation shows, that
\begin{equation}m=\min|f(x,a)|=0,\quad a>0,\quad  -\infty<x<+\infty.\label{3.5+}\end{equation}
Although the Karmarkar--Karp method gives the approximate  result $m=2$,
the precise  solution of the partition problem is $m=(9+7)-(6+5+5)=0$. Our analytic method gives (see \eqref{3.5+}) the same result $m=0$.
\section{Limit theorem}
\setcounter{equation}{0}
Let us consider the almost periodic functions of the form
\begin{equation}f(x)=\sum_{k=1}^{\infty}c_{k}e^{ix\lambda_k},\quad f_{n}(x)=\sum_{k=1}^{n}c_{k}e^{ix\lambda_k},
\label{4.1}\end{equation}
where
\begin{equation}\sum_{k=1}^{\infty}|c_{k}|<\infty.
\label{4.2}\end{equation}It follows from  \eqref{4.1} and \eqref{4.2} that
\begin{equation}|f(x)-f_{n}(x)|{\leq}\sum_{k=n+1}^{\infty}|c_{k}|{\to}0,\quad n{\to}\infty.
\label{4.3}\end{equation}
Using \eqref{4.3} we obtain
\begin{equation}\inf|f_{n}(x)|{\to}\inf|f(x)|,\quad \sup|f_{n}(x)|{\to}\sup|f(x)|,
\label{4.4}\end{equation}where $-\infty<x<+\infty,\quad n{\to}\infty.$
Taking into account \eqref{2.40} and \eqref{4.4} we derive our next assertion
\begin{Tm}\label{Theorem 4.1}Let the conditions \eqref{4.1} and \eqref{4.2} be fulfilled. Then,
\begin{equation}\min\Big|\sum_{k=1}^{n}(|c_k|\varepsilon_k)\Big|{\to}\inf|f(x)| \quad (-\infty<x<+\infty)\quad {\mathrm{for}} \quad n{\to}\infty.
\label{4.5}\end{equation} \end{Tm}
\begin{Rk}\label{Remark 4.2} In the proof of  Theorem \ref{Theorem 4.1}, we use the interrelations between the almost periodic functions theory and the partition theory. 
In our opinion, it would be  essentially more difficult to prove the existence of 
$$\lim_{n{\to}\infty}\min\Big|\sum_{k=1}^{n}(|c_k|\varepsilon_k)\Big|$$  
staying within the framework 
of the partition theory.\end{Rk}
\textbf{Acknowledgements}
I am deeply grateful to I. Tydniouk for the fruitful discussions, useful comments, and for the essential help in  charting and calculations. I am very grateful to A. Sakhnovich for
the  important help with the references, for a very productive discussion and for the improvements in the stile of this paper.
\appendix
\section{Some figures to examples}
\setcounter{equation}{0}

\begin{figure}[H]
  \centering
    \includegraphics[width=1.0\textwidth]{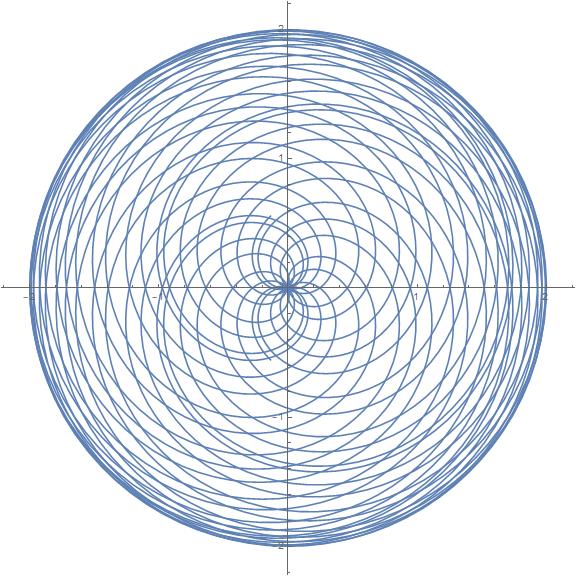}
		\label{fig:Example28}
  \caption{Example 2.8}
\end{figure}

\begin{figure}[H]
  \centering
    \includegraphics[width=1.0\textwidth]{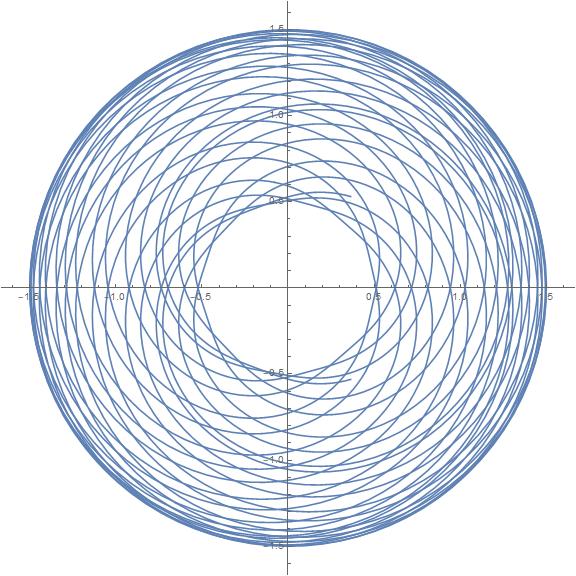}
		\label{fig:Example29}
  \caption{Example 2.9}

\end{figure}
\begin{figure}[H]
  \centering
    \includegraphics[width=1.0\textwidth]{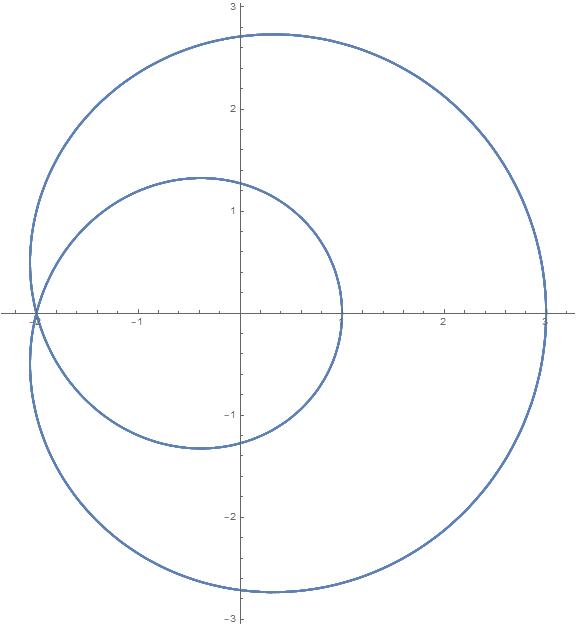}
		\label{fig:Example2.16}
  \caption{Example 2.16}
\end{figure}

\begin{figure}[H]
  \centering
    \includegraphics[width=1.0\textwidth]{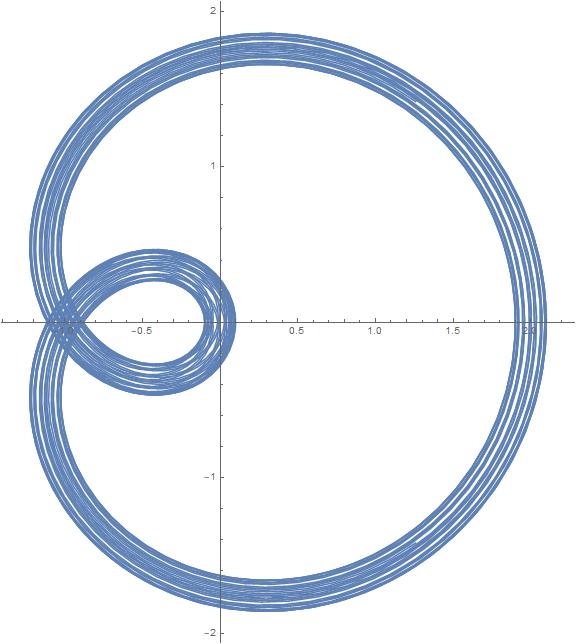}
		\label{fig:Example2.19}
  \caption{Example 2.19}
\end{figure}

\end{document}